\documentclass[a4paper,10pt]{article}

\usepackage[english]{babel}

\usepackage{lastpage}
\usepackage{amsmath,amsfonts,
  amssymb,amsthm,latexsym,dsfont,color,bbm,mathtools,mathrsfs}

\usepackage{geometry}
\RequirePackage{bera} 

\usepackage[utf8]{inputenc}
\usepackage[T1]{fontenc}
\usepackage[normalem]{ulem}

\usepackage{natbib}

\newcommand\unnumberedfootnote[1]{ %
  \let\temp=\thefootnote %
  \renewcommand{\thefootnote}{}%
  \footnote{#1}%
  \let\thefootnote=\temp%
  \addtocounter{footnote}{-1}}

\newcommand{\overbar}[1]{\mkern 1.5mu\overline{\mkern-1.5mu#1\mkern-1.5mu}\mkern 1.5mu}

\newcommand{\R}{\mathbbm{R}}

\newtheorem{theorem}{Theorem}
\newtheorem{proposition}{Proposition}[section]
\newtheorem{lemma}[proposition]{Lemma}
\newtheorem{corollary}[proposition]{Corollary}

\theoremstyle{definition}
\newtheorem{remark}[proposition]{Remark}
\numberwithin{equation}{section}
\numberwithin{equation}{section}

\begin{document}

\title{\LARGE Some large deviations in Kingman's coalescent}

\author{\sc Andrej Depperschmidt, Peter Pfaffelhuber \\
  \sc and Annika Scheuringer}

\vspace{2cm}
\date{\today}

\maketitle
\unnumberedfootnote{\emph{AMS 2010 subject classification.}  60F10, 60G09}
\unnumberedfootnote{\emph{Keywords and phrases.} Kingman coalescent,
  large deviations, uniform spacings}

\begin{abstract}
  Kingman's coalescent is a random tree that arises from classical
  population genetic models such as the Moran model. The individuals
  alive in these models correspond to the leaves in the tree and the
  following two laws of large numbers concerning the structure of the
  tree-top are well-known: (i) The (shortest) distance, denoted by
  $T_n$, from the tree-top to the level when there are $n$ lines in
  the tree satisfies $nT_n \xrightarrow{n\to\infty} 2$ almost surely;
  (ii) At time $T_n$, the population is naturally partitioned in
  exactly $n$ families where individuals belong to the same family if
  they have a common ancestor at time $T_n$ in the past. If $F_{i,n}$
  denotes the size of the $i$th family, then $n(F_{1,n}^2 + \cdots +
  F_{n,n}^2) \xrightarrow{n\to \infty}2$ almost surely. For both laws
  of large numbers we prove corresponding large deviations results.
  For (i), the rate of the large deviations is $n$ and we can give the
  rate function explicitly. For (ii), the rate is $n$ for downwards
  deviations and $\sqrt n$ for upwards deviations. For both cases we
  give the exact rate function.
\end{abstract}

\section{Introduction}
Kingman's coalescent is a random tree introduced by
\cite{Kingman1982a} as the genealogy arising in large population
genetic models. It has infinitely many leaves and is usually
constructed from leaves to the root as follows: given that there are
$k$ lines in the tree, after some exponential time with rate $\binom k
2$, two lines are chosen uniformly and merged to one line, leaving the
tree with $k-1$ lines. Due to the quadratic rate $\binom k2$ the tree
immediately comes down from infinitely to finitely many leaves
\citep{Donnelly1991}. Since the seminal paper by \cite{Pit1999} this
random tree has been generalized to other infinite trees arising in
population genetics models.

For the Kingman coalescent some laws of large numbers and central
limit theorems have been proved. They are nicely summarized in
\cite{Aldous1999}, Chapter~4.2; see also Proposition~\ref{P:11} below.
For $\varepsilon>0$ let $N_\varepsilon$ denote the number of lines
time $\varepsilon$ in the past. Then, since the Kingman coalescent
immediately comes down from infinity, $N_\varepsilon$ is finite.
Furthermore it is approximately $2/\varepsilon$. Equivalently, the
time $T_n$ it takes the coalescent to go from infinitely many lines to
$n$ lines is approximately $2/n$ for large $n$. Going to the fine
structure, at time $T_n$ the infinite population is decomposed in $n$
families (whose joint distribution is exchangeable) and every leaf in
the tree belongs to exactly one of the $n$ families whose frequencies
are denoted by $F_{1,n},\dots,F_{n,n}$. It is known that for large $n$
a randomly chosen $F_{i,n}$ is approximately exponentially distributed
with mean $1/n$. This translates into several laws of large numbers;
see e.g.\ (35) in \cite{Aldous1999}. In particular the probability of
picking (from the initial infinite population) two leaves that belong
to the same family, given by $F_{1,n}^2+ \dots+ F_{n,n}^2$, is
approximately $2/n$.

The main goal of the present paper is to study the corresponding large
deviations results. To the best of our knowledge, except for
\cite{AngelBerestyckiLimic2012}, cf.\ Remark~\ref{rem:angel}, results
in this direction are not present in the literature. We formulate our
results in the next section. Theorem~\ref{T1} gives a full large
deviation principle for the distributions of $nT_n$. The proof, given
in Section~\ref{sec:proofs1}, is an application of the G\"artner-Ellis
Theorem. As a byproduct, we derive a large deviation principle for the
distributions of $\varepsilon N_\varepsilon$ in
Corollary~\ref{cor:tNt}. Large deviations of $n(F_{1,n}^2+ \dots+
F_{n,n}^2)$ are considered in Theorem~\ref{T2} and exact rate
functions for downwards and upwards deviations are given. The proof is
given in Section~\ref{sec:proof-theorem-reft2}. For the upward
deviations we use a variant of Cramér's theorem for heavy-tailed
random variables; see e.g.~\cite{GantertRamananRembart2014}. For the
downward deviations we use a connection to self-normalized large
deviations; see \cite{Shao1997}. This connection was pointed out to us
by Alain Rouault and Nina Gantert. Since the rate function for
downward deviations is hard to treat analytically we provide in
Theorem~\ref{T3} a simple lower bound. The proof of that bound is
given in Section~\ref{sec:proof-theorem3}.

\section{Main results}
The Kingman coalescent can be seen as a discrete graph,
more precisely a discrete tree with infinitely many leaves.
Let $S_2, S_3,\dots$ be independent exponentially distributed
variables with mean $1$. Then the Kingman coalescent tree can be
constructed from the root to the leaves as follows.
\begin{enumerate}
\item Start the tree with two lines from the root.
\item For $k \ge 2$ the tree stays with $k$ lines for the amount of
  time $S_k/\binom k 2$. After that time one of the $k$ lines is
  randomly chosen.  This line splits in two so that the number of
  lines jumps from $k$ to $k+1$.
\item Stop upon reaching infinitely many lines, which happens after
  (the almost surely finite) time $T_1 \coloneqq \sum_{k=2}^\infty S_k/\binom
  k 2$.
\end{enumerate}
The random variable $T_1$ is the total tree height. Alternatively,
$T_1$ is the time to the most recent common ancestor (MRCA) of the
infinite population (of leaves). Counted from the top of the tree at
time $\varepsilon > 0$ a random number $N_\varepsilon$ of active lines
in the Kingman tree is present, i.e.\
\begin{align}
  \label{eq:Ne_Tn}
  N_\varepsilon \coloneqq \inf\{n: T_n<\varepsilon\} \qquad
  \text{for}\qquad T_n \coloneqq \sum_{k=n+1}^\infty \frac{S_k}{\binom
    k 2}.
\end{align}
At time $T_n$ every leaf belongs to one of $n$ disjoint families and
all members of each such family stem from the same line at time $T_n$.
Let us denote the frequencies of these families (which exist due to
exchangeability by deFinetti's Theorem) by
$F_{1,n},\dots,F_{n,n}$. The following results are well known (see
\cite{Aldous1999} for~\eqref{eq:lln-Tn} and \eqref{eq:lln-Ne} and
\cite{Evans2000} for \eqref{eq:lln-Fe}; proofs can also be found in
\cite{DepperschmidtGrevenPfaffelhuber2013}.)

\begin{proposition}[Laws of large numbers]
  \leavevmode \\
  Let \label{P:11} $(T_n)_{n=1,2,\dots}$,
  $(N_\varepsilon)_{\varepsilon>0}$ and
  $(F_{1,n},\dots,F_{n,n})_{n=1,2,\dots}$ be as above. Then
  \begin{align}
    \label{eq:lln-Tn}
    n T_n & \xrightarrow{n \to \infty} 2 \quad
    \text{almost surely,} \\
    \label{eq:lln-Ne}
    \varepsilon N_\varepsilon & \xrightarrow{\varepsilon \to 0} 2
    \quad \text{almost surely,} \\ \intertext{and}
    \label{eq:lln-Fe}
    n \sum_{k=1}^{n} F_{k,n}^2 & \xrightarrow{n\to\infty} 2 \quad
    \text{almost surely.}
  \end{align}
\end{proposition}

\begin{remark}[Interpretation of \eqref{eq:lln-Fe}]
  We \label{rem:interLLnFe} note that the left hand side of
  \eqref{eq:lln-Fe} has the interpretation of a {\it homozygosity by
    descent} in the following sense: when picking two leaves from the
  tree at time $0$, the probability that both share a common ancestor
  at time $T_n$ is $\sum_{k=1}^{n} F_{k,n}^2$. Then, the law of large
  number states that the homozygosity by descent at time $T_n$ is
  approximately $2/n$ for large $n$.
\end{remark}

\begin{figure}[t]\centering
  \includegraphics[width=6cm]{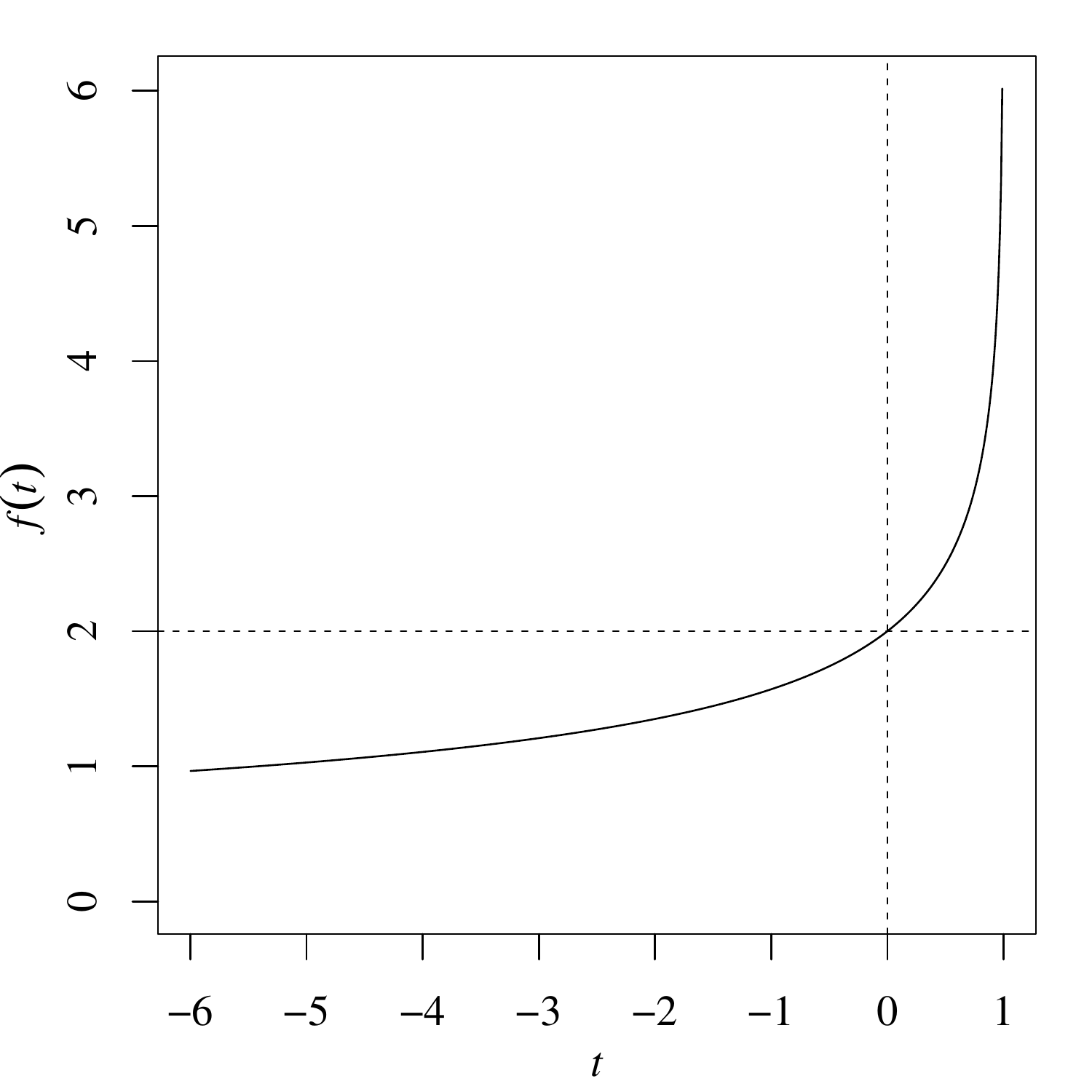} \hspace{1cm}
  \includegraphics[width=6cm]{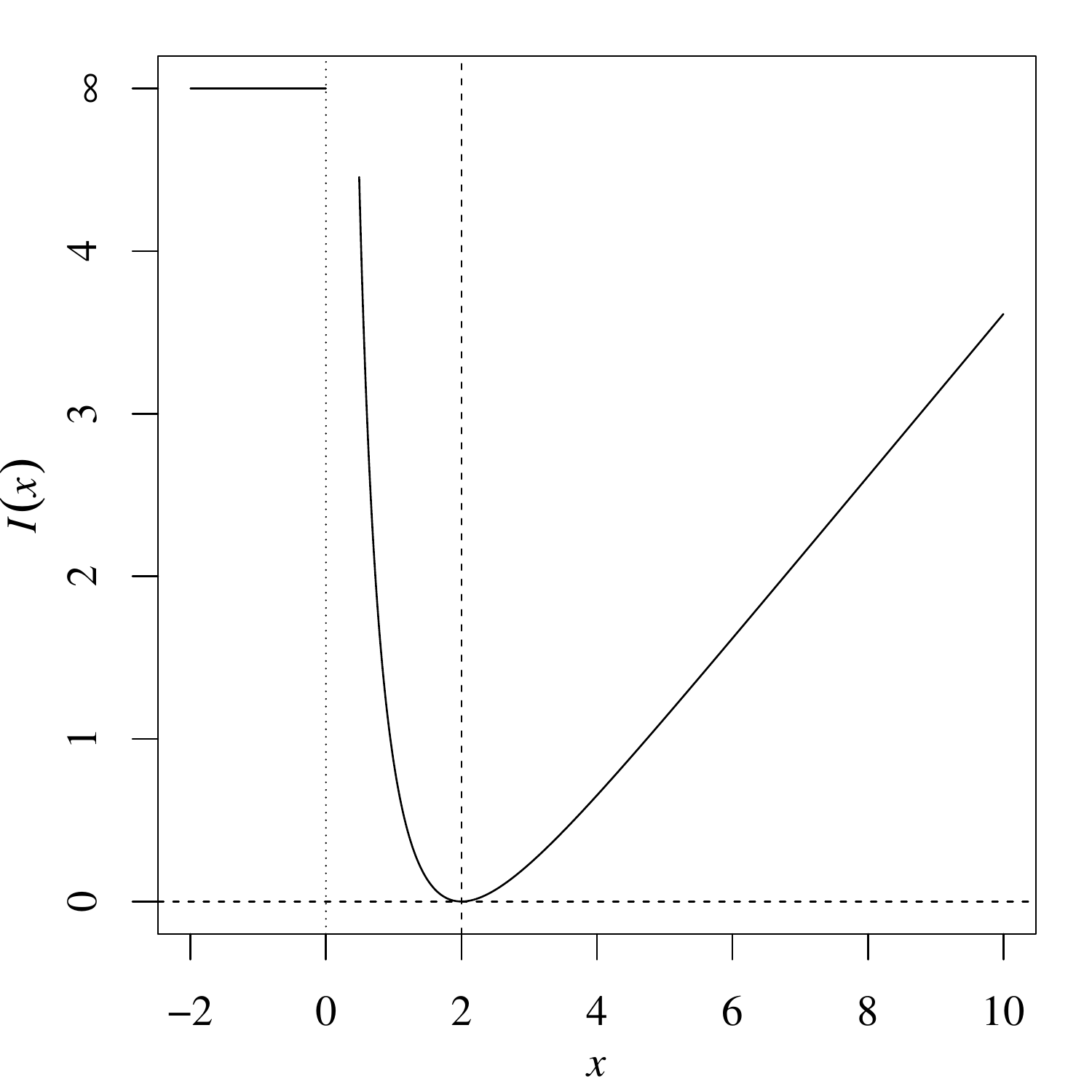}
  \caption{This figure displays the functions $f$ and $I$ from
    \eqref{eq:f} and \eqref{eq:I}, respectively.}
  \label{fig:T1}
\end{figure}

In the present paper we are interested in large deviations results
corresponding to the statements of Proposition~\ref{P:11}. We start
with large deviations connected with \eqref{eq:lln-Tn}. First we
introduce some notation. For $n=1,2,\dots$ let $\mu_n$ denote the
distribution of $nT_n$, i.e.\ $\mu_n( \, \cdot \, ) = \mathbf P ( n
T_n \in \, \cdot \, )$. Furthermore we denote by $\mathcal B(\mathbb
R)$ the Borel $\sigma$-algebra on $\R$ and for $\Gamma \in\mathcal
B(\mathbb R)$ we denote by $\Gamma^\circ$ the \emph{interior} and by
$\overbar \Gamma$ the \emph{closure} of $\Gamma$. For $x>0$, let
$t_x<1$ be the unique solution of the equation $x = f(t)$, where the
continuous and increasing function $f:(-\infty,1) \to (0,\infty)$ is
defined by (see Figure~\ref{fig:T1} for a plot)
\begin{align}\label{eq:f}
  f(t) \coloneqq \begin{cases} \displaystyle \frac{1}{\sqrt t} \log
    \frac{1+\sqrt t}{1-\sqrt t}  & : \; 0< t<1, \\[2ex]
    2 & : \; t=0, \\[2ex]
    \displaystyle\frac{2}{\sqrt{|t|}}\arctan\sqrt{|t|} & : \; t< 0.
  \end{cases}
\end{align}
The proof of the following theorem is given in
Section~\ref{sec:proof-theorem-reft1}.
\begin{theorem}[LDP for $(\mu_n)_{n=1,2,\dots}$] The sequence \label{T1}
  $(\mu_n)_{n=1,2,\dots}$ satisfies a large deviation principle with
  scale $n$ and good rate function $I$ given by
  \begin{align}\label{eq:I}
    I(x) \coloneqq
    \begin{cases}
      \displaystyle \frac{t_x}{2} x + \int_1^\infty \log\Big(1 -
      \frac{t_x}{y^2}\Big)\, dy & : \; x>0, \\[3ex]
      \infty & :\; x\leq 0.
    \end{cases}
  \end{align}
  In other words, for any $\Gamma \in\mathcal B(\mathbb R)$ we have
  \[ -\inf_{x\in \Gamma^\circ} I(x) \leq \liminf_{n\to\infty} \frac 1n
  \log \mu_n(\Gamma) \leq \limsup_{n\to\infty} \frac 1n \log \mu_n(\Gamma)
  \leq -\inf_{x\in \overbar \Gamma}I(x).\]
\end{theorem}

\begin{remark}[Interpretation]
  Both, the function $f$ from \eqref{eq:f} and $I$ from \eqref{eq:I}
  are plotted in Figure~\ref{fig:T1}. The minimum of the rate function
  is attained at $x=2$. This fact is clear from the law of large
  numbers, \eqref{eq:lln-Tn}. In addition, $I(x)=\infty$ for $x\leq 0$
  because $nT_n>0$ almost surely.

  \medskip
  Let us now have a closer look at the behaviour of $I(x)$ for $x$
  near $0$ and for large $x$. Since $2\arctan(t) \xrightarrow{t\to
    \infty} \pi$, we have that $\sqrt{|t_x|} x
  \xrightarrow{x\downarrow 0} \pi$, and hence, $t_x
  \stackrel{x\downarrow 0}\approx - \frac{\pi^2}{x^2}$. In this case,
  \begin{equation}
    \label{eq:I0}
    \begin{aligned}
      xI(x) & \stackrel{x\downarrow 0}\approx - \frac{\pi^2}{2} +
      x\int_1^{\infty} \log\Big(1 + \frac{\pi^2}{x^2 y^2}\Big) dy \\ &
      \stackrel{x\downarrow 0}\approx - \frac{\pi^2}{2} + \pi
      \int_{0}^\infty \log\Big(1+\frac{1}{z^2}\Big) dz =
      \frac{\pi^2}{2},
    \end{aligned}
  \end{equation}
  where the last equality follows from $\int_0^\infty \log(1+1/z^2)dz
  = \pi$. To understand the behaviour for large $x$, note that since
  \[ f(t) \stackrel{t\uparrow 1}\approx \log 2 - \log(1-\sqrt{t}),\]
  for $x \to \infty$ we have $2/(1-\sqrt{t_x}) \approx e^x$ and
  in particular $ t_x \approx (1-2e^{-x})^2 \approx 1$. It follows
  \[ \frac{I(x)}{x} \stackrel{x\to\infty} \approx \frac 12 +
  \frac{1}{x}\int_1^\infty \log(1-1/y^2) \,dy \stackrel{x\to\infty}
  \approx \frac 12.\]
\end{remark}

\begin{figure}[ht!]\centering
  \includegraphics[width=6cm]{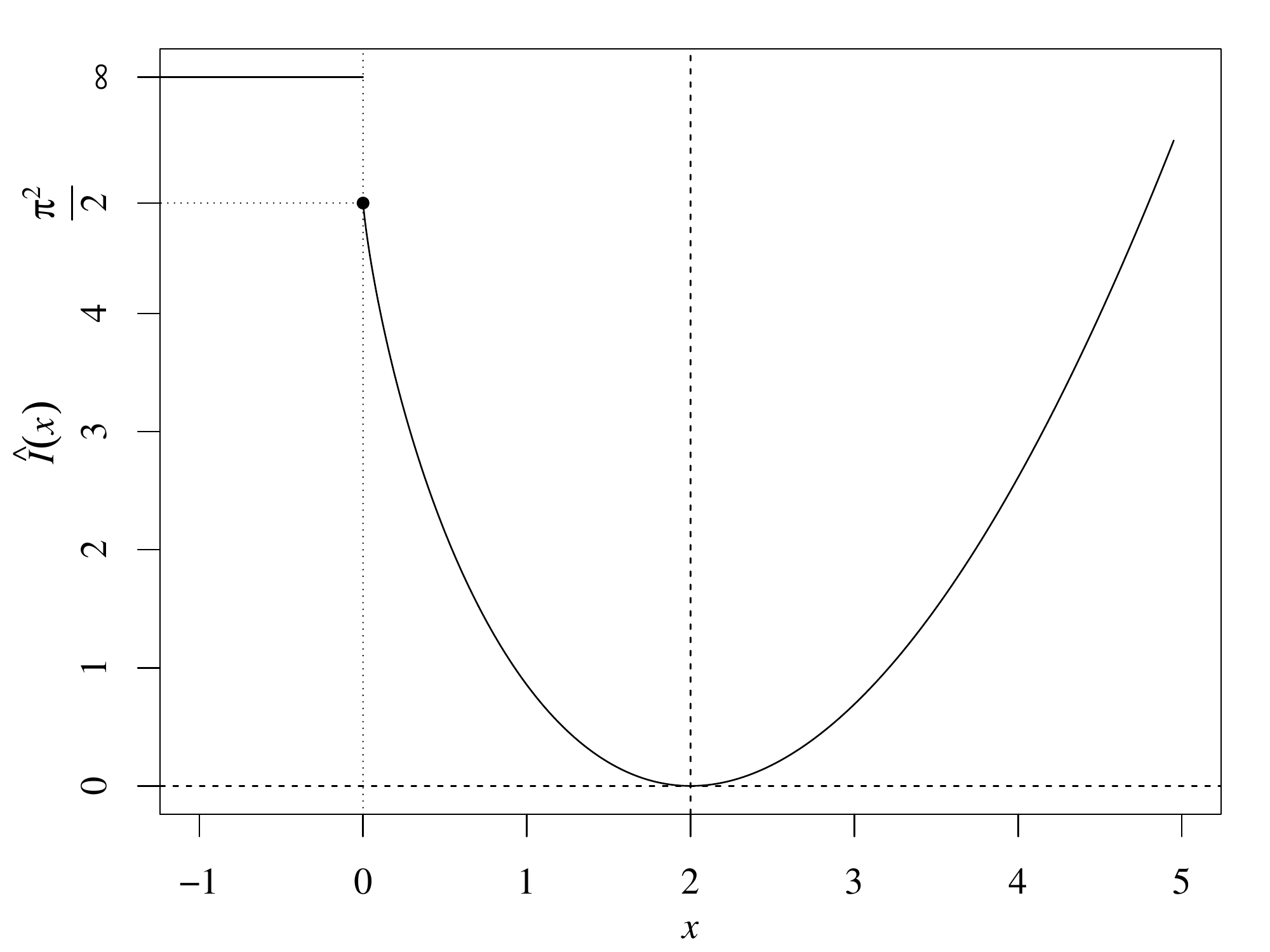} \hspace{1cm}
  \includegraphics[width=6cm]{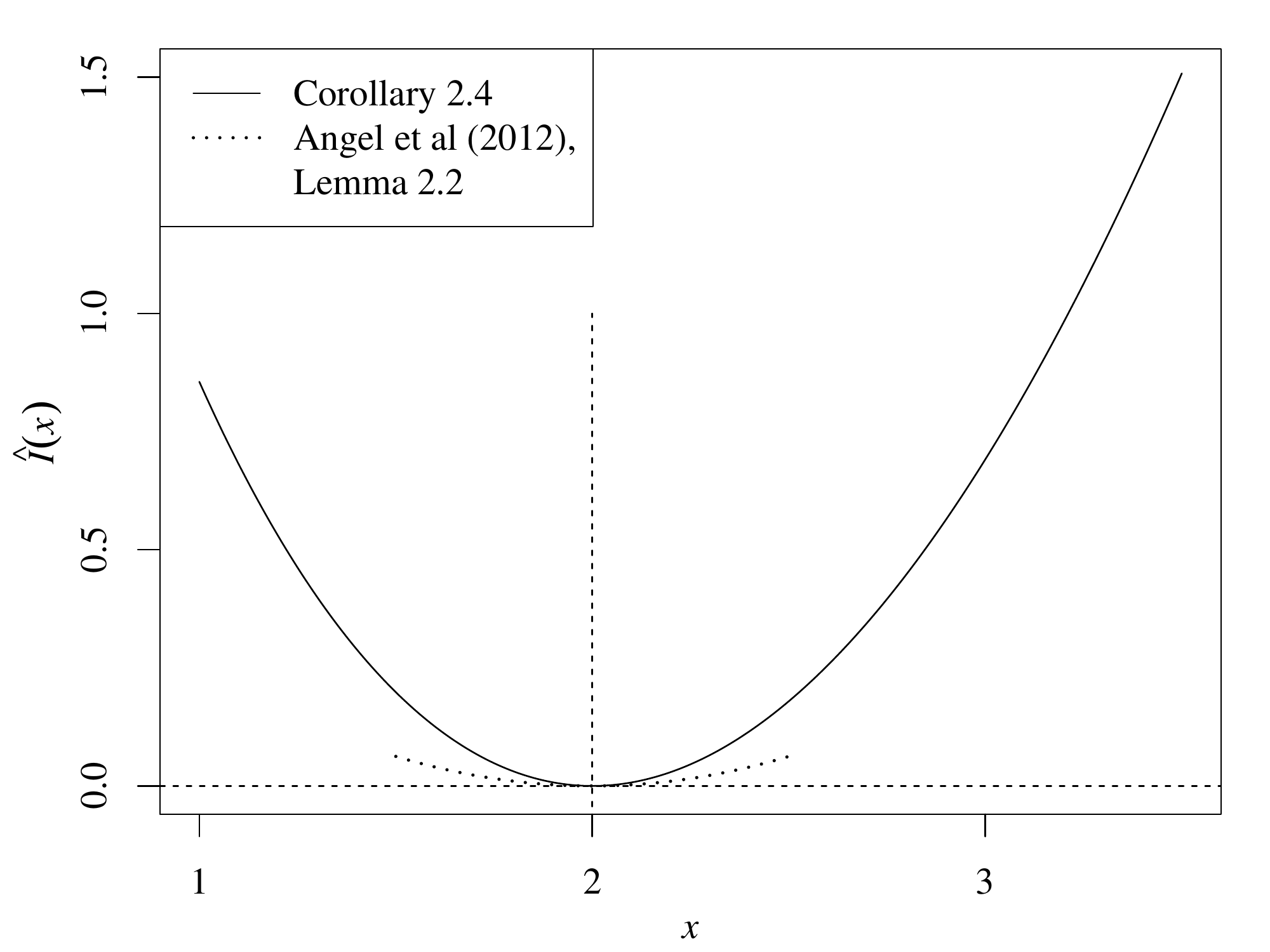}
  \caption{The figure on the left shows the rate function $\widehat I$
    from Corollary~\ref{cor:tNt}. The figure on the right is a
    comparison of $\widehat I$ with the lower bound obtained from
    \cite{AngelBerestyckiLimic2012}.}
  \label{fig:cor1}
\end{figure}

Note that \eqref{eq:lln-Tn} and \eqref{eq:lln-Ne} are equivalent.
Indeed, $\{T_n \ge \varepsilon\} = \{ N_\varepsilon \ge n\}$ (this
also holds with $\geq$ replaced by $\leq$) by construction, and $T_n
\downarrow 0$ as $n \to \infty$ and $N_\varepsilon \uparrow \infty$ as
$\varepsilon \to 0$. Hence, Theorem~\ref{T1} translates into a large
deviation principle for $\varepsilon N_\varepsilon$. In the following
we denote by $\nu_\varepsilon$ the distribution of $\varepsilon
N_\varepsilon$, i.e.\ $\nu_\varepsilon( \,\cdot\, ) = \mathbf
P(\varepsilon N_\varepsilon \in \, \cdot\,)$. The proof of the next
result is given in Section~\ref{sec:proof-coroll-refc}; see
Figure~\ref{fig:cor1} for a plot of the rate function $\widehat I$.

\begin{corollary}[LDP for $(\nu_\varepsilon)_{\varepsilon>0}$]
  For $ \varepsilon \downarrow 0$ the family
  $(\nu_\varepsilon)_{\varepsilon>0}$ \label{cor:tNt} satisfies a
  large deviation principle with scale $1/\varepsilon$ and good rate
  function $\widehat I$ given by
  \begin{align} \label{eq:Iprime}
    \widehat I(x) =
    \begin{cases}
      xI(x)  &: \, x>0,\\[1.5ex] \displaystyle
      \frac{\pi^2}{2} &: \,  x=0, \\[1.5ex]
      \infty &:\, x<0,
    \end{cases}
  \end{align}
  with $I$ from~\eqref{eq:I}. In particular, for $\Gamma \in\mathcal
  B(\mathbb R)$ we have
  \[ -\inf_{x\in \Gamma^\circ} \widehat I(x) \leq \liminf_{
    \varepsilon \to 0} \varepsilon  \log \nu_\varepsilon (\Gamma) \leq \limsup_{
    \varepsilon \to 0} \varepsilon  \log \nu_\varepsilon (\Gamma) \leq
  -\inf_{x\in \overbar \Gamma} \widehat I(x).\]
\end{corollary}

\begin{remark}[The full distribution of $N_\varepsilon$]
  The distributions $\nu_\varepsilon, \varepsilon>0$ (as well as
  $\mu_n, n=1,2,\dots$) have been described explicitely in the
  literature. \cite{Tavare1984}, Section 6, gives
  \[ \mathbf P(N_\varepsilon = n) = \sum_{k=n}^\infty e^{-\binom k 2 \varepsilon}
  \frac{(-1)^{k-n}(2k-1)\cdot n\cdots(n+k-2)}{n!(k-n)!}.\] In
  principle, this formula must also give the large deviations for the
  measures $\nu_\varepsilon$, but this does not seem straight-forward.
\end{remark}

\begin{remark}[The rate function $\widehat I$ and comparison
  with~\cite{AngelBerestyckiLimic2012}]
  Although \label{rem:angel} the main goal of~\cite{AngelBerestyckiLimic2012} was the
  analysis of spatial $\Lambda$-coalescents, they also provide some
  large deviations bounds on Kingman's coalescent. These bounds are
  mainly based on Markov inequality. Precisely, in Lemma 2.2 in
  \cite{AngelBerestyckiLimic2012} it is shown that for $0<x<\tfrac 12$
  \begin{align*}
    \mathbf P(|\varepsilon N_\varepsilon -2|>x) < e^{-\frac{x^2}{4\varepsilon}}
  \end{align*}
  and therefore
  \begin{align*}
     \limsup_{\varepsilon \to 0}  \varepsilon \log \mathbf
    P(| \varepsilon N_\varepsilon -2|>x) \leq - \frac{x^2}{4}.
  \end{align*}
  In the neighbourhood of $2$ the last inequality translates easily
  into a bound for the rate function $\widehat I$ from \eqref{eq:Iprime}; see
  Figure~\ref{fig:cor1}. Namely, for $x \in (1.5,2.5)$ we have
  \begin{align*}
    \widehat I(x) \geq \frac{(x-2)^2}{4}.
  \end{align*}
\end{remark}

\medskip
Next, we state some large deviations results connected
to~\eqref{eq:lln-Fe}. For
\begin{align*}
  W_n \coloneqq n\sum_{k=1}^{n} F_{k,n}^2
\end{align*}
we know from~\eqref{eq:lln-Fe} that $W_n\xrightarrow{n\to\infty} 2$
holds almost surely. The proof of this result is based on the
well-known fact (see e.g.\ Section~5 in \cite{Kingman1982a}) that the
distribution of $W_n$ can be derived using uniform order statistics:
Let $U_1,\dots,U_{n-1}$ be independent and uniformly distributed
on $[0,1]$, and $0<U_{(1)} < \cdots < U_{(n-1)} < 1$ be their order
statistics.  Additionally, let $R_1,\dots,R_{n}$ be independent
exponentially distributed random variables with mean $1$. Then,
\begin{equation}
  \label{eq:762}
  \begin{aligned}
    (F_{1,n},\dots,F_{n,n}) & \stackrel{d} = \bigl(U_{(1)},
    U_{(2)}-U_{(1)},\dots,U_{(n-1)} - U_{(n-2)}, 1-U_{(n-1)}\bigr)
    \\ & \stackrel d = \Big( \frac{R_1}{\sum_{j=1}^n
      R_j},\dots,\frac{R_n}{\sum_{j=1}^n R_j}\Big).
  \end{aligned}
\end{equation}
Here the second equality in distribution is one of the well known
representations of uniform spacings; see e.g.\ Section~4.1 in
\cite{Pyke1965}. It follows
\begin{align}\label{eq:763}
  W_n \stackrel d = n \frac{\sum_{k=1}^n R_k^2}{\Big(\sum_{j=1}^n
    R_j\Big)^2} = \frac{\frac 1n \sum_{k=1}^n R_k^2}{\Big(\frac 1n
    \sum_{j=1}^n R_j\Big)^2}.
\end{align}
We will use this representation to obtain large deviations results for
$W_n$. In particular we show that upwards large deviations of $W_n$
are on the scale $\sqrt n$ while downwards large deviations are on the
scale $n$. The proof is given in
Section~\ref{sec:proof-theorem-reft2}.

\begin{theorem}[Large deviations of $W_n$]
  For \label{T2} each $x \ge 2$, we have
  \begin{align}\label{eq:ldWn1}
    \lim_{n\to\infty} \frac{1}{\sqrt{n}} \log \mathbf{P}\left( W_n
      \geq x \right) = - \sqrt{x-2}.
  \end{align}
  Furthermore $\mathbf P(W_n<1)=0$ and for each $ 1 < x < 2$, we have
  \begin{align}
    \label{eq:ldWn-down}
    \lim_{n\to\infty} \frac{1}{n} \log \mathbf{P}\left( W_n \leq x
    \right) = - \widetilde I (x).
  \end{align}
  The function $\widetilde I(x)$ is positive for $1 <x <2$ and is
  given by
  \begin{align}
    \label{eq:tildeI}
  \widetilde I (x) \coloneqq \sup_{c \ge 0} \inf_{t \ge 0} M(x,c,t).
  \end{align}
  Here $M: (1,2) \times[0,\infty) \times [0,\infty) \to
  \R$ is a function of the form $M\coloneqq M_1 + M_2+M_3$ with
  \begin{align*}
    M_1(x,c,t) & \coloneqq \frac12\log(2\pi) +
    \frac14\log x-\frac12\log t\\
    M_2(x,c,t) & \coloneqq \frac{(tc-1)^2 x -t^2c^2}{2t\sqrt{x}}\\
    M_3(x,c,t) & \coloneqq  \log
    \Phi\left(\frac{(tc-1)x^{1/4}}{\sqrt t}\right)
  \end{align*}
  where $\Phi$ denotes the distribution function of the one
  dimensional standard Gaussian distribution.
\end{theorem}

Though the rate function in~\eqref{eq:ldWn-down} is exact it is hard
to treat analytically. For this reason we provide in Theorem~\ref{T3}
a much simpler lower bound for downwards large deviations of $W_n$.
For the proof we use the following lemma which provides another
representation of $W_n$ in terms of exponential random variables (see
Section~\ref{sec:proofs2} for proofs).
\begin{lemma}[Representation of $W_n$]
  Let \label{l:Wn} $R_1,\dots,R_n$ be independent exponentially
  distributed random variables with mean $1$. Then,
  \begin{align}
    \label{eq:Wn-repr}
     W_n \stackrel d = \frac{\frac 1n
    \Big(2\sum_{l=1}^n \sum_{k=1}^l \frac{R_kR_l}{l} - \sum_{k=1}^n
    \frac{R_k^2}{k}\Big)}{\Big(\frac 1n \sum_{k=1}^n R_k\Big)^2}.
  \end{align}
\end{lemma}

\begin{theorem}[Lower bound on downwards large deviations of $W_n$]
  For \label{T3} $1 < x < 2$ we have
  \begin{align}
    \label{eq:ldWn2}
         \liminf_{n\to\infty} \frac{1}{n} \log \mathbf{P}\left( W_n
           \leq x \right) \geq    1 - \frac{1}{\sqrt{x-1}}.
  \end{align}
\end{theorem}

\begin{remark}[Rationale and use of the representation in
  Lemma~\ref{l:Wn}]
  The main point in the proof of Lemma~\ref{l:Wn} is that $W_n$ does
  not depend on the order of the $R_k$ and hence we can as well order
  them according to their size.

  Let us briefly explain how we will use \eqref{eq:Wn-repr} in the
  proof of in \eqref{eq:ldWn2}. Since $W_n$ is minimal if $R_1=\cdots
  = R_n$ (whence $W_n=1$), we have to look for possibilities that all
  $R_k$'s are of about the same size in order to obtain a large
  deviations result for $W_n$. Let $R_{(1)}, \dots, R_{(n)}$ denote
  the above exponential random variables ordered in increasing order,
  i.e.\ $R_{(i)}$ is the $i$th smallest value. Using ``competing
  exponential clocks'' arguments (see also the proof of the lemma) one
  can see that $R_{(i)} - R_{(i-1)}$ is exponentially distributed with
  mean $1/(n-i+1)$. Hence, one way of obtaining similar values for all
  $R_k$'s arises if $R_{(1)}$ is particularly large, which then leads
  to a large deviations result for $W_n$.
\end{remark}

\begin{figure}[ht!]\centering
  \includegraphics[width=6cm]{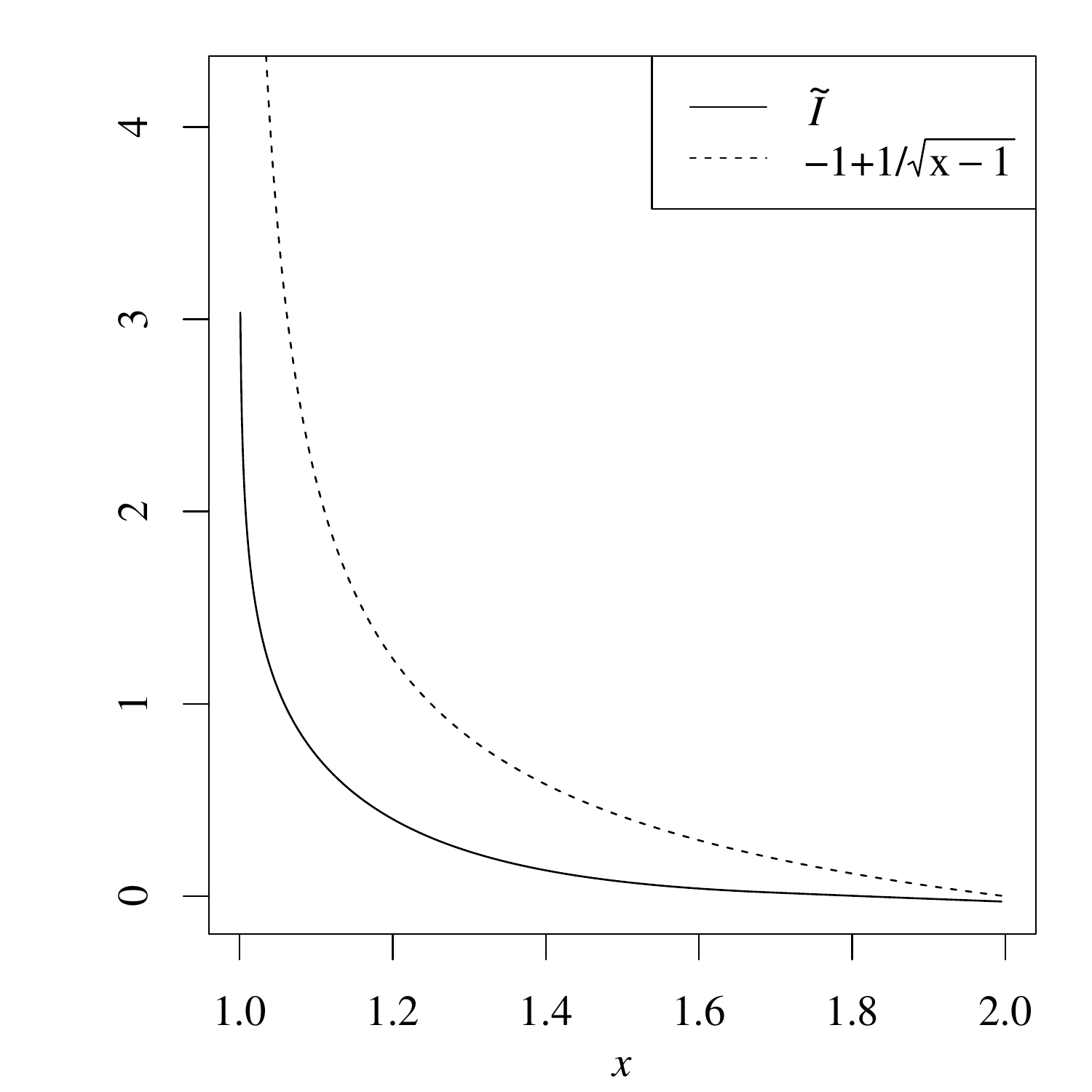}
  \caption{Numerical comparison of the exact rate function for
    downwards large deviations of $W_n$ from \eqref{eq:ldWn-down} in
    Theorem~\ref{T2} and the lower bound from
   \eqref{eq:ldWn2} in Theorem~\ref{T3}.}
  \label{fig:T2a}
\end{figure}

\begin{remark}[Interpretation, uniform spacings and the
  Poisson-Dirichlet distribution]\leavevmode
  1.\/ Let us give some heuristics about the rates arising in
  Theorem~\ref{T2}. For \eqref{eq:ldWn1}, we have to ask ourselves
  about the easiest way $W_n$ becomes too large. From \eqref{eq:762},
  we see that this is the case if one of the $R_k$'s is too large,
  making this kind of deviations a local property in the sense that
  only a single of the $R_k$'s has to show some untypical behavior.
  This is different when looking at \eqref{eq:ldWn-down}, i.e.\ too
  small values of $W_n$. First, observe that $W_n$ is small only if
  all (or many) families have about equal sizes (extreme case $F_{1,n}
  = \cdots = F_{n,n} = \frac 1n$ gives the minimal value $W_n=1$).
  Hence, such downward deviations require to study a global property
  of the random variable $W_n$, which is significantly harder. For the
  proof of \eqref{eq:ldWn-down} we will interpret $W_n$ as a
  self-normalised sum and use from \cite{Shao1997} a result on large
  deviations result for such sums.

  \medskip \noindent
  2.\/ From \eqref{eq:762}, we see that in fact $W_n$ is a function of
  uniform order statistics, which, for instance, have been studied in
  detail (although no large deviations results were given) in
  \cite{Pyke1965}. Hence, Theorem~\ref{T2} may as well be interpreted
  as a large deviations result for uniform order statistics.

  \medskip \noindent 3.\/ As stated in Remark~\ref{rem:interLLnFe},
  $W_n /n$ can be interpreted as homozygosity at time $T_n$. Using a
  Poisson process along the tree with intensity $\theta/2$, we can ask
  for the probability of picking two leaves from the tree which are
  not separated by a Poisson mark, denoted by {\it homozygosity in
    state}, abbreviated by $H_\theta/\theta$. This quantity is closely
  related to the Poisson-Dirichlet distribution and some large
  deviations (in the limit of large $\theta$) were derived in
  \cite{DawsonFeng2006}. It is shown there in Theorem~5.1 that
  $H_\theta/\theta \xrightarrow{\theta\to\infty}0$ and that
    \begin{align*}
      \mathbf P(H_\theta > \theta x) = e^{-\theta (I(\sqrt x) + o(1))}
    \end{align*}
    for $I(x) = -\log(1-x)$. However, a large deviation principle for
    the quantity $H_\theta$ (noting that $H_\theta
    \xrightarrow{\theta\to\infty}1$), which corresponds to the results
    from Theorem~\ref{T2}, could not be obtained by
    \cite{DawsonFeng2006}. At least, it was shown that its scale
    cannot be larger than $\sqrt\theta$.
\end{remark}

\section[Proof of Theorem~\ref{T1} and Corollary~\ref{cor:tNt}]{Proof
  of Theorem~\ref{T1} and Corollary~\ref{cor:tNt}}
\label{sec:proofs1}
\subsection{Proof of Theorem~\ref{T1}}
\label{sec:proof-theorem-reft1}

The proof of Theorem~\ref{T1} is an application of the G\"artner-Ellis
theorem; see for instance Section~2.3 in \cite{DemboZeitouni2010}.

Let $\Lambda_n(t) \coloneqq \log \mathbf E[e^{t nT_n}]$ and $\mu_n( \,
\cdot \, ) = \mathbf P ( n T_n \in \, \cdot \, )$. To show that the
sequence $\mu_1, \mu_2,\dots$ satisfies a large deviation principle
with scale $n$ and a good rate function we need to check the following
three conditions.
\begin{enumerate}
\item[GE1] $\Lambda(t) \coloneqq \lim_{n\to\infty} \frac 1n
  \Lambda_n(nt)$ exists for all $t$ as a limit in $\overline{\mathbb
    R} = \mathbb R \cup\{\pm \infty\}$. Furthermore $t \to \Lambda
  (t)$ is lower-semicontinuous, $0 \in \mathcal
  D_\Lambda^\circ$, where $\mathcal D_{\Lambda} \coloneqq
  \{t: \Lambda(t)<\infty\}$.

\item[GE2] $\Lambda$ is differentiable on $\mathcal
  D_\Lambda^\circ$.
\item[GE3] $\Lambda$ is \emph{steep}, i.e.\ $\Lambda'(t_n)
  \xrightarrow{n\to\infty}\infty$ whenever $t_1,t_2,\dots\in \mathcal
  D_{\Lambda}^\circ$ and $t_n\xrightarrow{n\to\infty} t\in \partial
  \mathcal D_{\Lambda}$.
\end{enumerate}
Then the good rate function is given by
\begin{align}
  \label{eq:gRF}
  x\mapsto I(x) = \sup_{t\in\mathbb R}(tx - \Lambda(t)).
\end{align}

We proceed in three steps. First, we compute $\Lambda(t) \coloneqq
\lim_{n\to\infty} \frac1n\Lambda_n(nt)$. Second, we check the further
assumptions of the G\"artner-Ellis theorem and obtain $I$ as the
Fenchel-Legendre transform of $\Lambda$. In the third step, for the
rate function $I$ from \eqref{eq:gRF} we obtain its simplified form
given in Theorem~\ref{T1}.

\medskip
\noindent {\bf Step 1. The limit of $\boldsymbol{\frac1n\Lambda_n(nt)}$:}
We will show that
\begin{align}
  \label{eq:Lambda-limit}
  \Lambda(t) = \lim_{n\to\infty} \frac{1}{n} \Lambda_n (nt) =
  \begin{cases}
    - \int_1^\infty \log\left(1- \frac{2t}{x^2}\right) \, dx & : t
    \le \frac12,   \\
    \infty & : t > \frac12.
  \end{cases}
\end{align}
For this, recall from \eqref{eq:Ne_Tn} that $T_n = \sum_{k=n+1}^\infty
S_k/\binom k 2$ where $S_k/\binom k2$ is exponentially
distributed with rate $\binom{k}{2}$ as well as independent of
$S_\ell$ for all $\ell \neq k $. Furthermore recall that the moment
generating function of an exponentially distributed random variable
$R$ with rate $\lambda>0$ is given by
\begin{align}\label{expmomallg}
  \mathbf E\left[e^{tR}\right] =
  \begin{cases}
    \frac{\lambda}{\lambda-t}, &\text{if $t<\lambda$},\\ \infty,
    &\text{if $t\geq \lambda$.}
  \end{cases}
\end{align}
Hence, for each $n\in\mathbb N$ and $t\in \mathbb R$ we obtain by the
monotone convergence theorem
\begin{equation}\label{expmom}
  \varphi_n (nt) \coloneqq \mathbf{E}\left[e^{t n^2 T_n}\right]
  = \mathbf{E}\left[e^{t n^2 \sum_{k=n+1}^\infty S_k/\binom k2}\right]
  = \prod_{k=n+1}^\infty \mathbf{E}\left[e^{tn^2 S_k/\binom k 2}\right].
\end{equation}
We have to consider two cases $t > \frac12$ and $t \le \frac12$
separately. First suppose that $t>\frac{1}{2}$. Then there exists $n_0
\in \mathbb{N}$ so that for all $n\geq n_0$ we have
\begin{align*}
  2t \geq \frac{n+1}{n}, \quad \text{i.e.}\quad
  \frac{tn^2}{\binom{n+1} 2} = 2t \frac{n}{n+1}>1.
\end{align*}
Consequently, using \eqref{expmomallg}, we obtain $\mathbf E[e^{tn^2 T_{n}}]
= \infty$ for each $n\geq n_0$. Hence, $\varphi_n (nt) = \infty$ and
$\Lambda_n (nt)= \log \varphi_n (nt) =\infty$ for $n$ large enough.
Thus, we have
\begin{align*}
  \Lambda(t) =\lim_{n\to\infty}\frac{1}{n} \Lambda_n (nt) = \infty
  \quad \text{for all } \;  t > \frac{1}{2}.
\end{align*}
Now suppose that $t \leq \frac{1}{2}$. For $n\in \mathbb N$ and $k\geq n+1$
we have $ \binom{k}{2} \geq tn^2$. Furthermore using \eqref{expmom}
and \eqref{expmomallg} we can write
\begin{align*}
  \varphi_n (nt) = \prod_{k=n+1}^\infty \frac{1}{1-\frac{2 t
      n^2}{k(k-1)}}.
\end{align*}
Using this we can rewrite $ \frac{1}{n}\Lambda_n (nt)$ for $t\leq
\tfrac 12 $ as
\begin{equation*}
  \begin{aligned}
    \frac{1}{n}\Lambda_n (nt)
    &= \frac{1}{n}\log\left( \prod_{k=n+1}^\infty \frac{1}{1- \frac{2 t n^2}{k(k-1)}}\right)
    = - \frac{1}{n}\sum_{k=n}^\infty \log
    \left(1-\frac{2t}{\frac{k}{n}\frac{k+1}{n} }\right) \\
    &=- \frac{1}{n}\sum_{x\in \left\{1,1+\frac{1}{n},
        1+\frac{2}{n},\dots \right\}} \log \left( 1-
      \frac{2t}{x\left(x+\frac{1}{n}\right)} \right)
    = - \int_1^\infty \log \left(1- \frac{2t}{\frac{\lfloor xn
          \rfloor}{n} \frac{\lfloor xn \rfloor+1}{n}} \right) \,dx
  \end{aligned}
\end{equation*}
and by the dominated convergence theorem we obtain
\begin{align*}
    \frac{1}{n}\Lambda_n (nt) \xrightarrow{n\to\infty} - \int_1^\infty \log \left(1-
      \frac{2t}{x^2} \right) \,dx.
\end{align*}
Hence, GE1 is shown with $\Lambda$ as in \eqref{eq:Lambda-limit}.
Moreover, we have $\mathcal D_\Lambda = (-\infty, \tfrac 12]$,
$\Lambda(\tfrac 12) = \int_1^\infty \log(1-\tfrac 1{x^2},dx = \pi$ and
$\Lambda$ is lower-semi-continuous.

\medskip
\noindent {\bf Step 2. Further assumptions of the G\"artner-Ellis
  theorem:}
We proceed by checking the assumptions GE2 and GE3. For
differentiability of $\Lambda$ for $t<\frac{1}{2}$ consider for
$-\infty < r < 0 < s < \frac{1}{2}$ the function
\begin{align*}
  f:(1,\infty) \times (r,s) &\to \mathbb{R} \\
  (x,t)\hspace{0.65cm} &\mapsto -\log\left(1-\frac{2t}{x^2}\right).
\end{align*}
We have $\int_1^\infty \left| f(x,t) \right| \, dx < \infty$ for
$t \in (r,s)$ and the derivative
\begin{align*}
  \frac{d}{dt} f(x,t) = \frac{2}{x^2-2t}
\end{align*}
exists for each $x \in (1,\infty)$ and is continuous in $t$.  Hence,
we can interchange differentiation and integration and obtain
\begin{align*}
  \Lambda'(t) = \int_1^\infty \frac{2}{x^2-2t} \, dx.
\end{align*}
Furthermore, for a sequence $t_1, t_2,\dots$ with $t_n \uparrow \frac12$
we obtain
\begin{align*}
  \lim_{n \to\infty }\Lambda'(t_n) & = \lim_{n \to\infty}
  \int_1^\infty \frac{2}{x^2-2t_n} \, dx \\
  &= \lim_{n \to \infty} \frac{1}{\sqrt{2t_n}} \left(
    \log\left(1+\sqrt{2t_n}\right) -
    \log\left(1-\sqrt{2t_n}\right)\right) =\infty,
\end{align*}
i.e.\ condition GE3 is also satisfied.

\medskip
\noindent {\bf Step 3. Properties of $\boldsymbol{I}$:} Applying the
G\"artner-Ellis theorem reveals that the sequence of distributions of
$nT_n, n=1,2,\dots$ satisfies a large deviation principle with good
rate function
\begin{align*}
  I(x) & = \sup_{t \leq \tfrac 12} \left[ tx + \int_1^\infty
    \log\left(1-\frac{2t}{y^2}\right) \,\emph dy \right]= \sup_{t \leq
    1} \left[ \frac t2 x + \int_1^\infty
    \log\left(1-\frac{t}{y^2}\right) \,\emph dy \right].
\end{align*}
In order to compute that supremum, we write for $t\geq 0$
\begin{align*}
  2 \frac{\partial}{\partial t} \left[ \frac t2 x + \int_1^\infty
    \log\left(1-\frac{t}{y^2}\right) \, dy\right] & = x -
  2\int_1^\infty \frac{1}{y^2 - t}dy \\ & = x +
  \frac{1}{\sqrt{t}}\int_1^\infty \frac{1}{y+\sqrt{t}} -
  \frac{1}{y-\sqrt{t}} dy \\ & = x - \frac{1}{\sqrt{t}} \log
  \frac{1+\sqrt{t}}{1-\sqrt{t}}
\end{align*}
while for $t\leq 0$
\begin{align*}
  2 \frac{\partial}{\partial t} \left[ \frac t2 x + \int_1^\infty
    \log\left(1-\frac{t}{y^2}\right) \, dy\right] & = x -
  2\int_1^\infty \frac{1}{y^2 + |t|}dy \\ & = x -
  \frac{2}{\sqrt{|t|}}\arctan\sqrt{|t|}.
\end{align*}
It is easy to see that the second derivative is negative throughout,
such that the supremum is attained at $t_x$ given by the solution of
$f(t_x)=x$ for $f$ as in~\eqref{eq:f}. Finally we note that for $t\in
[0,1)$ the range of $t\mapsto \frac{1}{\sqrt{t}} \log
\frac{1+\sqrt{t}}{1-\sqrt{t}}$ is $[2,\infty)$ and for $t\in
(-\infty,0]$ the range of $t\mapsto
\frac{2}{\sqrt{|t|}}\arctan\sqrt{|t|}$ is $(0,2]$. Hence, the scale
function $I$ is of the form given in~\eqref{eq:I}.

\subsection{Proof of Corollary~\ref{cor:tNt}}
\label{sec:proof-coroll-refc}
The proof is based on the fact that $\{T_n \ge \varepsilon\} = \{
N_\varepsilon \ge n\}$. Thus, for $x\geq 2$ we have
\begin{align*}
  x I(x) & = x \lim_{n\to\infty} \frac 1n \log \mathbf P(nT_n \geq x)
  = \lim_{n\to\infty} \frac xn \log \mathbf P(\frac xn N_{x/n} \geq x)
  = \lim_{\varepsilon \to 0} \varepsilon \log \mathbf P(\varepsilon
  N_{\varepsilon} \geq x)
\end{align*}
and for $0<x\leq 2$
\begin{align*}
  x I(x) & = x \lim_{n\to\infty} \frac 1n \log \mathbf P(nT_n \leq
  x) = \lim_{n\to\infty} \frac xn \log \mathbf P(\frac xn N_{x/n}
  \leq x) = \lim_{\varepsilon\to 0} \varepsilon \log \mathbf
  P(\varepsilon N_{\varepsilon} \leq x).
\end{align*}
The value $\widehat I(0)$ follows from~\eqref{eq:I0}. Since the rate
function $I$ attains its minimum at $x=2$, is decreasing below and
increasing above $2$, the result follows.

\section{Proof of Lemma~\ref{l:Wn}, Theorem~\ref{T2} and
  Theorem~\ref{T3} }
\label{sec:proofs2}

\subsection{Proof of Lemma~\ref{l:Wn}}
\label{sec:proof-lemma-refl:wn}
When looking at \eqref{eq:763}, note that $W_n$ does not depend on the
order of the $R_k$'s. Therefore, it is possible to order them
according to their size. Precisely, let $0 < R_{(1)} < \cdots <
R_{(n)}$ be their order statistics. Then it is well-known that
\begin{align*}
  (R_{(n)},R_{(n-1)},\dots,R_{(1)}) \stackrel d = \left(\sum_{k=1}^n
    \frac{R_k}{k}, \sum_{k=2}^n \frac{R_k}{k},\dots,
    \frac{R_n}{n}\right), \quad \text{i.e.} \quad R_{(n-k+1)}
  \stackrel d = \sum_{i=k}^n \frac{R_i}{i}.
\end{align*}
Indeed, the smallest of $n$ independent exponentially distributed mean
$1$ random variables is exponentially distributed with mean $\tfrac1n$ (as
does $\tfrac{R_n}n$), and the second smallest then has the same distribution
as $\tfrac{R_n}{n} + \tfrac{R_{n-1}}{n-1}$ etc. Now, we
obtain \eqref{eq:Wn-repr} as follows
\begin{align*}
  W_n & \stackrel d = \frac{\frac 1n \sum_{k=1}^n
    R_{(n-k+1)}^2}{\Big(\frac 1n \sum_{j=1}^n R_{(n-j+1)}\Big)^2}
  \stackrel d = \frac{\frac 1n \sum_{k=1}^n \Big(\sum_{i=k}^n
    \frac{R_i}{i}\Big)^2}{\Big(\frac 1n \sum_{j=1}^n \sum_{i=j}^n
    \frac{R_i}{i}\Big)^2} \\ & = \frac{\frac 1n \Big( 2 \sum_{k=1}^n
    \sum_{i=k}^n \sum_{j=i}^n \frac{R_i R_j}{ij} - \sum_{k=1}^n
    \sum_{i=k}^n \frac{R_i^2}{i^2}\Big)}{\Big( \frac 1n \sum_{i=1}^n
    \sum_{j=1}^i \frac{R_i}{i}\Big)^2} \\ & = \frac{\frac 1n
    \Big(2\sum_{j=1}^n \sum_{i=1}^j \sum_{k=1}^i \frac{R_iR_j}{ij} -
    \sum_{i=1}^n \sum_{k=1}^i \frac{R_i^2}{i^2}\Big)}{\Big(\frac 1n
    \sum_{i=1}^n R_i\Big)^2} \\ & = \frac{\frac 1n \Big(2\sum_{j=1}^n
    \sum_{i=1}^j \frac{R_iR_j}{j} - \sum_{i=1}^n
    \frac{R_i^2}{i}\Big)}{\Big(\frac 1n \sum_{i=1}^n R_i\Big)^2}.
\end{align*}

\subsection{Proof of Theorem~\ref{T2}}
\label{sec:proof-theorem-reft2}
We start by proving~\eqref{eq:ldWn1}. Let $x \ge 2$ and let $R_1,R_2,
\dots$ be independent exponential random variables with mean $1$. In
what follows we set
\begin{align}
  \label{eq:R-abbr}
  X_n \coloneqq \frac1n\sum_{k=1}^n R_k \quad \text{and} \quad Z_n
  \coloneqq \frac{1}{n} \sum_{k=1}^n R_k^2.
\end{align}
According to \eqref{eq:763}, it suffices to show that
\begin{align}\label{darstellungexp}
  \lim_{n\to\infty} \frac{1}{\sqrt{n}} \log \mathbf{P}\left(
    \frac{Z_n}{X_n^2} \geq x\right) = - \sqrt{x-2}.
\end{align}
To this end we will show that for all $0< \varepsilon <1$,
\begin{align}\label{fall1toshow}
  \limsup_{n\to\infty} \frac{1}{\sqrt{n}} \log \mathbf{P}\left(
    \frac{Z_n}{X_n^2} \geq \frac{x}{1-\varepsilon} \right) \leq -
  \sqrt{x-2}
\end{align}
as well as
\begin{align}\label{fall2toshow}
  \liminf_{n\to\infty} \frac{1}{\sqrt{n}} \log \mathbf{P}\left(
    \frac{Z_n}{X_n^2} \geq x \right) \geq - \sqrt{(1+\varepsilon)x -2}
\end{align}
and obtain \eqref{darstellungexp} by letting $\varepsilon\to 0$.
For~\eqref{fall1toshow} we have
\begin{equation}\label{fall1}
  \mathbf{P}\left(\frac{Z_n}{X_n^2} \geq
    \frac{x}{1-\varepsilon} \right) \leq
  \mathbf{P}\left(Z_n \geq x\right) + \mathbf{P}\left(X_n \leq
    \sqrt{1-\varepsilon}\, \right).
\end{equation}
We consider the two terms on the right hand side of the last display
separately and start with the first one. Observe that $\mathbf
E[e^{\lambda R_1^2}] = \infty$ for $\lambda>0$, $\mathbf{E}[R_1^2] =2$
and $ \mathbf{P}\left(R_1^2 \geq t\right) = e^{-\sqrt{t}}$ for $t\geq
0$. We use a variant of Cram\'er's theorem for heavy-tailed
random variables from \cite{GantertRamananRembart2014}. In particular,
we refer to the statement around equation (1.2) there (the assumption
there is fulfilled with $X_1$ replaced by $R_1^2$ and $r=\frac{1}{2}$,
$m=2$ and $c=1$). We obtain
\begin{align}\label{p1fall1}
  \mathbf{P}\left(Z_n \geq x\right) =
  e^{-\sqrt{n} \left(\sqrt{x-2} + o(1)\right)} \quad \text{as
    $n\to\infty$}.
\end{align}
For the second term on the right hand side of \eqref{fall1} by the
(classical) Cramér theorem we obtain
\begin{align}\label{p2fall1}
  \mathbf{P}\left(X_n \leq \sqrt{1-\varepsilon} \right) = e^{-n
    I_{\exp}(\sqrt{1-\varepsilon})(1+o(1))}, \quad \text{as
    $n\to\infty$,}
\end{align}
where
\begin{align}
  \label{eq:Iexp}
  I_{\exp}(y) \coloneqq y - 1 -\log(y)
\end{align}
is the Fenchel-Legendre transform of the function $t\mapsto \log
\mathbf E[e^{\lambda R_1}]$. Now, using \eqref{fall1}, \eqref{p1fall1}
and \eqref{p2fall1} we obtain
\begin{multline*}
  \limsup_{n\to\infty} \frac{1}{\sqrt{n}} \log \mathbf{P}\left(
    \frac{Z_n}{X_n^2} \geq \frac{x}{1-\varepsilon}
  \right) \\ \leq \limsup_{n\to\infty} \frac{1}{\sqrt{n}} \log
  \Big(e^{-\sqrt{n}(\sqrt{x-2} + o(1))} +
  e^{-nI_{\exp}(\sqrt{1-\varepsilon})(1+o(1))}\Big) = - \sqrt{x-2},
\end{multline*}
which shows \eqref{fall1toshow}. For the proof of~\eqref{fall2toshow}
we write
\begin{equation}\label{fall2}
  \begin{split}
    \mathbf{P}\left( \frac{Z_n}{X_n^2} \geq x
    \right) &\geq \mathbf{P}\left( \frac{Z_n}{X_n^2} \geq x,
     X_n^2 \leq 1+\varepsilon\right) \\
    &\geq \mathbf{P}\left( Z_n \geq x(1+\varepsilon), X_n^2
      \leq 1+\varepsilon\right)
    \\
    &\geq\mathbf{P}\left( Z_n \geq x(1+\varepsilon) \right) -
    \mathbf{P}\left(X_n \geq \sqrt{1+\varepsilon}\right).
  \end{split}
\end{equation}
Again we consider both terms in the last line separately.
For the first term, as in \eqref{p1fall1} we obtain
\begin{align}\label{p1fall2}
  \mathbf{P}\left(Z_n \geq x(1+\varepsilon)\right) = e^{-\sqrt{n}
    \left(\sqrt{x(1+\varepsilon) - 2} + o(1)\right)}, \quad \text{as
    $n\to\infty$}.
\end{align}
For the second term, we use the same argument as for \eqref{p2fall1}
and get
\begin{align}\label{p2fall2}
  \mathbf{P}\left(X_n \geq \sqrt{1+\varepsilon}\right) \leq e^{-n
    I_{\exp}(\sqrt{1+\varepsilon})(1+o(1))}, \quad \text{as
    $n\to\infty$}.
\end{align}
Combining~\eqref{p1fall2} and~\eqref{p2fall2} with~\eqref{fall2} now
gives~\eqref{fall2toshow} which proves~\eqref{eq:ldWn1}.

\medskip
Since the minimum of $W_n$ is $1$ (when $F_{k,n} = 1/n$ for all $k$)
the assertion $\mathbf{P}(W_n<1)=0$ is clear. It remains to
prove \eqref{eq:ldWn-down}, show that the rate function is of the
form~\eqref{eq:tildeI} and justify the positivity of $\widetilde I(x)$
for $x \in (1,2)$.

For $x\in (1,2)$ using \eqref{eq:763} we obtain
\begin{align}
  \label{eq:selfnorm}
  \mathbf P (W_n \le x) = \mathbf P \left(\frac{\sum_{j=1}^n
      R_j}{\sqrt{n} \sqrt{\sum_{k=1}^n R_k^2}} \ge
    \frac1{\sqrt{x}}\right).
\end{align}
Furthermore, for $x \in (1,2)$ we have $1/\sqrt{x} >1/\sqrt{2}=\mathbf
E[R_1]/\sqrt{\mathbf E[R_1^2]}$. Thus, we can use Theorem~1.1 from
\cite{Shao1997} and obtain
\begin{align}
  \label{eq:shao_to_ld}
  \mathbf P (W_n \le x)^{1/n} = \sup_{c\ge 0} \inf_{t \ge 0} \mathbf
  E\left[ \exp\Bigl(t(c R_1 - \frac1{2\sqrt
      x}(R_1^2+c^2))\Bigr)\right].
\end{align}
Now we have
\begin{align*}
  E\left[ \exp\Bigl(t(c R_1 - \frac1{2\sqrt
      x}(R_1^2+c^2))\Bigr)\right] & = \int_0^\infty \exp\Bigl(-y + t(c
  y - \frac1{2\sqrt x}(y^2+c^2))\Bigr) \, dy
\\ \intertext{and elementary integration yields}
& = \frac{\sqrt{2\pi}x^{1/4}}{\sqrt{t}}
    \exp\left(\frac{(tc-1)^2x-t^2c^2}{2t\sqrt{x}}\right)
    \Phi\left(\frac{(tc-1)x^{1/4}}{\sqrt t}\right),
\end{align*}
where $\Phi$ denotes the distribution function of the one dimensional
standard Gaussian distribution. Taking $\log$ of the last term we
obtain \eqref{eq:tildeI}.

Now we fix $x \in (1,2)$ and show that $\widetilde I(x)$ is positive.
In the sequel we write \[h(r,c)\coloneqq
cr-\tfrac1{2\sqrt{x}}(r^2+c^2).\]
We have
\begin{align*}
\inf_{t \ge 0} \mathbf E & \left[ \exp\Bigl(t h(R_1,c) \Bigr)\right]
\\ & \ge  \mathbf
  E\left[ \inf_{t \ge 0}\exp\Bigl(t h(R_1,c) \Bigr)\right] \\
  & = \mathbf E\left[\mathbbm 1_{\{h(R_1,c) <0\}} \inf_{t \ge
      0}\exp\Bigl(t h(R_1,c) \Bigr)\right] +
  \mathbf E\left[ \mathbbm 1_{\{h(R_1,c) \ge 0\}}  \inf_{t \ge
      0}\exp\Bigl(t h(R_1,c) \Bigr)\right] \\
  &= \mathbf P \bigl(h(R_1,c) \ge 0\bigr).
\end{align*}
The function $r \mapsto h(r,c)$ is non-negative on the interval
$[r_1,r_2]$ where $r_{1/2} =
r_{1/2}(c) \coloneqq c(\sqrt{x}\pm\sqrt{x-1})$ are the zeros of the
function. It follows
\begin{align*}
  \mathbf E\left[ \inf_{t \ge 0}\exp\Bigl(t h(R_1,c) \Bigr)\right] =
  \mathbf P \left( r_1 \le R_1 \le r_2 \right) =
  e^{-c(\sqrt{x}-\sqrt{x-1})} - e^{-c(\sqrt{x}+\sqrt{x-1})}.
\end{align*}
Finally, by elementary calculation we obtain
\begin{align*}
   \sup_{c \ge 0}\bigl( e^{-c(\sqrt{x}-\sqrt{x-1})} -
   e^{-c(\sqrt{x}+\sqrt{x-1})}\bigr)= \frac{2\sqrt{x-1}}{\sqrt x +
     \sqrt{x-1}} \biggl(\frac{\sqrt{x}-\sqrt{x-1}}{\sqrt{x}+\sqrt{x-1}}
   \biggr)^{\frac{\sqrt{x}-\sqrt{x-1}}{2\sqrt{x-1}}}.
\end{align*}
This expression (and therefore also $\widetilde I(x)$) is positive for
$x \in (1,2)$. Thus, the proof of Theorem~\ref{T2} is concluded.

\subsection{Proof of Theorem~\ref{T3}}
\label{sec:proof-theorem3}

We prove the inequality \eqref{eq:ldWn2} using Lemma~\ref{l:Wn}. Let
$1< x<2$ and set $y= \frac{1}{\sqrt{x-1}}-1$. For $\varepsilon>0$ we
have
\begin{align*}
  \frac{1}{n} & \log \mathbf{P}\big(W_n \leq
  x+\varepsilon\big)\\
  &=\frac{1}{n} \log \mathbf{P}\bigg(n \frac{2 \sum_{l=1}^n
    \sum_{k=1}^l \frac{R_k R_l}{l} - \sum_{k=1}^n
    \frac{R_k^2}{k}}{\left(\sum_{k=1}^n
      R_k\right)^2} \leq x+\varepsilon\bigg)  \\
  &\geq \frac{1}{n} \log \mathbf{P}\bigg(n \frac{2 \sum_{l=1}^n
    \sum_{k=1}^l \frac{R_k R_l}{l} - \sum_{k=1}^n
    \frac{R_k^2}{k}}{\left(\sum_{k=1}^n R_k\right)^2} \leq
  x+\varepsilon,
  R_n > ny\bigg) \\
  & = \frac{1}{n} \log \Bigg\{\mathbf{P}\big(R_n > ny\big)\,
  \mathbf{P}\bigg( n\frac{2 \sum_{l=1}^n \sum_{k=1}^l \frac{R_k
      R_l}{l} - \sum_{k=1}^n \frac{R_k^2}{k}}{\left(\sum_{k=1}^n
      R_k\right)^2} \leq x+ \varepsilon \,\Bigg|\, R_n >
  ny\bigg)\Bigg\}.
\end{align*}
Now $\frac1n \log \mathbf{P}\big(R_n > ny\big)=-y$, and conditioning in
the second factor in the curly braces can be removed by using the fact
that conditioned on $R_n>ny$ the exponential random variable $R_n$ has
the same distribution as $ny+R_n$. After some elementary calculations
we see that the last line of the above display equals
\begin{align*}
  -y + \frac{1}{n} \log \mathbf{P}\Bigg(\frac{\frac{1}{n} \Big( 2
    \sum_{l=1}^n \sum_{k=1}^l \frac{R_k R_l}{l} - \sum_{k=1}^n
    \frac{R_k^2}{k}\Big) + 2y \frac{1}{n} \sum_{k=1}^n R_k +
    y^2}{\left(\frac{1}{n} \sum_{k=1}^n R_k\right)^2 + 2y
    \frac{1}{n} \sum_{k=1}^n R_k + y^2}\leq
  x+\varepsilon\Bigg).
\end{align*}
From the strong law of large numbers and \eqref{eq:lln-Fe} with
Lemma~\ref{l:Wn} we know that
\begin{align*}
\frac{1}{n} \sum_{k=1}^n R_k \xrightarrow{n\to\infty} 1, \quad
\text{and} \quad \frac{1}{n} \Big(2 \sum_{l=1}^n \sum_{k=1}^l
\frac{R_k R_l}{l} - \sum_{k=1}^n \frac{R_k^2}{k}\Big) &
\xrightarrow[]{n\to\infty} 2 \quad \text{almost surely.}
\end{align*}
It follows that almost surely
\begin{align*}
  \frac{\frac{1}{n} \Big( 2 \sum_{l=1}^n \sum_{k=1}^l \frac{R_k
      R_l}{l} - \sum_{k=1}^n \frac{R_k^2}{k}\Big) + 2y \frac{1}{n}
    \sum_{k=1}^n R_k + y^2}{\left(\frac{1}{n} \sum_{k=1}^n
      R_k\right)^2 + 2y \frac{1}{n} \sum_{k=1}^n R_k + y^2}
  \xrightarrow[]{n\to\infty} \frac{2+2y+y^2}{1+2y+y^2} = x.
\end{align*}
Thus,
\begin{align*}
  \liminf_{n\to\infty} \frac{1}{n} \log \mathbf{P}\Big( W_n \leq
  2-x+\varepsilon\Big) \geq -y = 1- \frac{1}{\sqrt{x-1}}.
\end{align*}
The rest follows by letting $\varepsilon\downarrow 0$.

\subsubsection{Acknowledgments}
We thank Shui Feng for pointing out connections to
\cite{DawsonFeng2006} and Nina Gantert and Alain Rouault for pointing
out the reference \cite{Shao1997} and fruitful email discussion that
led to the exact rate function in \eqref{eq:ldWn-down}. This research
was supported by the DFG through grants Pf-672/6-1 to AD and PP.


\end{document}